\documentclass[a4paper,12pt]{article}

\usepackage{a4wide}  
\usepackage{multicol,amssymb,stmaryrd}

\begin{document}

\vspace{1cm}
\begin{center}
{\Large \bf Sheaf Representations
for Algebraic Systems}\\[5mm]{\large A personal historical account\\[5mm]
Klaus Keimel}
\end{center}

\vspace{2cm}
I am reporting here how I was involved in sheaf
representations of various structures between roughly 1965 and 1975. During
that period I was working in close collaboration with Karl Heinrich
Hofmann who has greatly influenced my work.  \\

{\bf Background}\\

Pierre Shapira says that \emph{sheaves on topological spaces were
introduced by Jean Leray as a tool to deduce global properties from
local ones}. Leray's worked on these ideas as a prisoner of war
between 1940 and 1945. He published his definition of a sheaf in 1946
in the Comptes Rendus Acad. Sc. Paris 232, p.1366. Leray's definition
was modified by Lazard in Expos\'e 14, S\'eminaire Cartan
1950-51. Sheaves were introduced as a methodological
tool in algebraic topology, in particular in cohomology theory. The
extension of certain results 
on cohomology from compact spaces to locally compact spaces was only
possible by admitting `variable coefficient groups depending
continuously on the points of the space' as opposed to the
usual constant cohomology groups. A second source
for sheaves was the theory of analytic functions (sheaves of germs
of analytic functions). In 1953 they were
finally introduced in algebraic topology. These developments
culminated in the work of Grothendieck. A major source is:
\begin{verse}
 A. Grothendieck and J. Dieudonn\'e, \'El\'ements de g\'eom\'etrie
 alg\'ebrique. I. Le langage des sch\'emas, Inst. Hautes Etudes
 Sci. Publ. Math. No. 4 (1960), 228 pages.   
\end{verse}
The early history of the development of the notion of a sheaf with
an extensive bibliography one can find in: 
\begin{verse}
J. W. Gray, Fragments of the history of sheaf theory, in:
M. P. Fourman, C.~J. Mulvey, and D.~S. Scott (eds.), Applications
of Sheaves, Lecture Notes in Mathematics 753 (1979), 1 -- 79, 
Springer-Verlag.  
\end{verse}
Gray's contribution shows that the development in the line
of `representations of algebras by sheaves' (which I am referring
to in the sequel) is a minor one compared with the main research
topics. \\ 

{\bf Sheaf representations}\\

It seems to me that Grothendieck's \emph{affine scheme} associated with a
commutative ring with unit was the motivating example for the
developments that I was involved in. Whereas the developments
mentioned in the previous section
were driven by methodological needs of classical theories,
the point of view of the community 
that emerged around 1965 was that of `Representation of algebraic
structures by sections in sheaves'. A typical theorem reads as
follows
\begin{verse}
\emph{Let $A$ be an algebraic structure. There is a sheaf $\mathcal S$ of
algebraic structures $A_x, x\in X,$ of the same type over a
topological space $X$ such that 
$A$ is isomorphic to the algebra $\Gamma(\mathcal S)$ of all global
sections (alternatively: all global sections with compact support)}.
\end{verse}   
The hope is that the stalks $A_x$ are of a simpler nature than
$A$ itself and that one will gain a better understanding of the
structure of $A$ by 
means of the structure of the stalks $A_x$ and the structure of the
base space $X$. But in this community, there was no particular problem that
urged the use of sheaves.   \\

{\bf Sheaf representations of rings, the Zariski topology}\\

At the beginning the efforts were concentrated on rings, mostly
commutative rings with or without units. Later -- with less success --
non-commutative rings were attacked.

The base space $X$ was generally obtained as the space of all prime ideals with
the Zariski topology or as a subspace thereof.\\ 

Personally I entered the world of sheaf representations after
reading:
\begin{verse}
J. Dauns and K. H. Hofmann, The representation of biregular rings
by sheaves, Math. Z. 91 (1966), 103-123. \\
\end{verse}
A ring (not necessarily commutative or unital) is \emph{biregular},
if every ideal is generated by a central idempotent. The ideals of
such a ring form a generalized Boolean algebra (a distributive lattice
in which every principal ideal is a Boolean algebra). The prime ideals
with the Zariski topology form a locally Boolean space. The associated
sheaf is a sheaf of simple rings with identity, and the given ring is
isomorphic to the ring of all sections with compact support. In the
commutative case the stalks are fields, \\

The Memoir:
\begin{verse}
R. S. Pierce, Modules over commutative regular rings, Memoirs AMS 70
(1967), 112pp.
\end{verse}
is a step further in the same direction and was greatly influential at
the time.\\

I entered the scene by observing that the methods developed in the
previous papers could also be applied to lattice-ordered rings: 
\begin{verse}
K. Keimel, Repr\'esentation d'anneaux r\'eticul\'es dans des faisceaux,
C. R. Acad. Sci. Paris Sér. A-B 266 (1968), A124-A127. MR 37 6224.\\
K. Keimel, Anneaux r\'eticul\'es quasi-r\'eguliers et hyperarchim\'ediens,
C. R. Acad. Sci. Paris Sér. A-B 266 (1968), A524-A525. MR 37 6225.
\end{verse}
I had indeed learned about ordered algebraic structures during my year
at Tulane University in 1964/65 from Lectures delivered by A. H. Clifford. And
when I came to Paris in 1965, I joined the group around
M.-L. Dubreil-Jacotin (comprising Alain Bigard and Samuel Wolfenstein)
working on lattice-ordered groups. \\

{\bf Lattice-ordered groups and rings}\\

I then observed that the ring structure was
not needed for the sheaf representations and I extended my work to
commutative lattice ordered groups (which can be viewed as lattice ordered rings
with the multiplication being identically zero). This led to: 
\begin{verse}
K. Keimel, Repr\'esentation de groupes et d'anneaux r\'eticul\'es par des
sections dans des faisceaux, Th\`ese d'\'Etat, Paris 1970.\\
\end{verse}
On my homepage
(http://www.mathematik.tu-darmstadt.de/~keimel/publications.html) 
one can find the contents of this thesis, partly without
proofs. Chapter I deals with prime elements of lattices and the
Zariski (= Stonean) topology on them. Chapter II deals with the
Stone-\v Cech-Isbell 
compactification.  Chapter III finally presents various sheaf
representations. They are formulated for lattice-ordered \emph{annelo\"ides}:
These are lattice ordered `rings' with a non-commutative addition. In
this way one also includes non-abelian lattice
ordered groups as lattice-ordered annelo\"ides with trivial
multiplication. This latter aspect has been omitted in:  
\begin{verse}
K. Keimel, The representation of lattice-ordered groups and rings by
sections in sheaves, Lecture Notes in Math. 248, Springer-Verlag,
Berlin and New York, 1971, pages 1--98.
\end{verse}
This article is an improved version of Chapter III of the thesis
above. It contains: \\
- A general representation theorem, \\
- Representation over the direct factor space,\\
-  Representation over the space of minimal irreducible
$\ell$-ideals,\\
- Representation of quasi-regular and hyperarchimedean f-rings and
$\ell$-groups,\\  
- Representation over the maximal $\ell$-ideal space, Dominated
$\ell$-ideals. \\ 

Some of the results of the article just indicated are reproduced in:
\begin{verse}
A. Bigard, K. Keimel, S. Wolfenstein, Groupes et Anneaux R\'eticul\'es
Lecture Notes in Mathematics, 608 (1977), Springer Verlag, xiv+334
pages.
\end{verse}

One should not forget the following contribution along the same lines:
\begin{verse}
J. Dauns, Representation of L-groups and F-rings, Pacific
J. Math. 31 (1969), 629-654: MR 41 130.
\end{verse}
It should be noted that a research group in Roumania had already
begun to work on sheaf representation of distributive lattices:  
\begin{verse}
A. Brezuleanu and R. Diaconescu, Sur la duale de la cat\'egorie des
treillis, Rev. Roumaine Math. Pures Appl. 14 (1969), 311-323. MR 41 
116.\\
A. Brezuleanu, Sur les sch\'emas de treillis, Rev. Roumaine
Math. Pures Appl. 14 (1969), 949 - 954.
\end{verse}

{\bf Sheaf representations for universal algebra}\\

The extension of the sheaf representation methods from rings to
lattice-ordered structures 
lead me to generalize further and to publish a representation theorem
for universal algebras, in particular, semigroups, based on a
generalized Boolean lattice of commuting factor congruences:
\begin{verse}
K. Keimel, Darstellung von Halbgruppen und universellen Algebren durch
Schnitte in 
Garben; biregul\"are Halbgruppen, Math. Nachrichten 45 (1970), 81 - 96.
\end{verse}
As this paper was published in German, it has  not been cited very often. 
Shortly later, but independently, Comer published a paper that is more
easily accessible and written in a better style based on the same ideas:  
\begin{verse}
S. D. Comer, Representation by algebras of sections over Boolean
spaces, Pacific J. Math. 38 (1971), 29 - 38.
\end{verse}
In the following years these ideas were developed further in numerous
papers as for example in:
\begin{verse}
K. Keimel and H. Werner, Stone duality for varieties generated by
quasi-primal algebras, 
Memoirs Amer. Math. Soc. 148 (1974), pp. 59-85. 
\end{verse}
In this context one should mention:
\begin{verse}
B.~A. Davey, Sheaf spaces and sheaves of universal algebras,
Math. Z. 134 (1973), 275-290.
\end{verse}
A student of mine wrote an interesting paper on sheaf representations
of universal algebras in general (the title is misleading as it seems
to point to classical algebras):
\begin{verse}
A. Wolf, Sheaf representations of arithmetical algebras, Memoirs
Amer. Math. Soc. 148 (1974), pp. 85-97. 
\end{verse}
In this paper a general setting is described that allows 
sheaf representations of universal 
algebras: \emph{The essential ingredient for a good sheaf representation is
the existence of a distributive lattice of pairwise commuting
congruences.}
It seems to be the first time that \emph{the Chinese Remainder
Theorem is shown to be responsible for the patch property, which
allows us to glue sections and which is
essential for proving that the given algebra is represented by the
algebra of ALL global sections (of compact support) of the associated
sheaf.} But note that appropriate formulations of the Chinese Remainder
Theorem for universal algebras were already available in the late
1950-ies; see for example the book on Universal Algebra by G. Gr\"atzer.\\

Independently, quite similar observations have been made by Swamy, who
reproves the appropriate Chinese Remainder Theorem: 
\begin{verse}
I. M. Swamy, Representation of universal algebras by sheaves,
Proc.~Amer. Math. Soc. 45 (1974),  55-58.
\end{verse}
When I came to Darmstadt in 1971, there was a strong group working in
Universal Algebra and Lattice Theory around R. Wille. After my
arrival, I introduced sheaf representations there. H. Werner and
A. Wolf began to work together with me. S. Burris and
B. Davey were visiting Darmstadt and contributed to the spreading of 
these techniques around to the community of Universal Algebra.\\

{\bf A Survey}\\

The following survey contains the developments on sheaf
representations described 
above in much more detail
up to the year 1972 with a rich list of references: 
\begin{verse}
K. H. Hofmann, Representation of algebras by continuous sections,
Bull. Amer. math. Soc. 78, Number 3
(1972), 291-373. \\
---, Some bibliographical remark on "Representation of algebras by
continuous sections", memoirs Amer. Math. Soc. 148 (1974), 177-182. 
\end{verse}

{\bf What do I mean by sheaf representation?}\\

By a \emph{sheaf} I here mean a local homeomorphism $p\colon E\to X$
of topological spaces. The \emph{stalks} are the $E_x=p^{-1}(x);x\in
X$. The n-fold fibered product is the sheaf $p^{(n)}\colon E^{(n)}\to
X$, where $E^{(n)} = \{(a_1,\dots,a_n)\mid p(a_1)=\dots=p(a_n)\}$ with
the subspace topology from $E^n$. 

If $\Omega$ is a finitary signature, a sheaf of
$\Omega$-algebras is a sheaf $p\colon E\to X$ together with continuous
maps $\omega\colon E^{(n)}\to E$ which respect fibers, where
$\omega\in\Omega$ is of arity $n$. This means in particular that the
stalks $E_x$ are $\Omega$-algebras.\\ 

For an $\Omega$-algebra $A$, a \emph{sheaf representation} is a an
algebra isomorphism of $A$ onto the algebra $\Gamma(p)$ of all global
(continuous) sections of a sheaf $p\colon E\to X$ of $\Omega$-algebras.\\

In order to find a sheaf representation of an $\Omega$-algebra $A$,
one may proceed roughly as follows. One chooses a set $X$ of congruence
relations $c_x$ of $A$. One forms the disjoint union $E$ of all the
quotient algebras $A_x=A/c_x, x\in X,$ together with the obvious projection
$p\colon E\to X$. For every $a\in A$ let $\widehat a\colon X\to E$
be defined by $\widehat a(x)=a\mbox{ mod }c_x$. One then chooses a
topology on $X$. 

1) If for all $a,b \in A$, the 'equalizer' $e(a,b)=\{x\in X\mid
(a,b)\in c_x\}$ is 
open, the we lift this topology up to $E$ by declaring the sets
$\widehat a(e(b,c))$ to be open. We then have a sheaf of
$\Omega$-algebras $p\colon E\to X$ where the stalks are the quotients
$A/c_x$. The map $a\mapsto\widehat a$ is an algebra homomorphisms  
of $A$ into the algebra $\Gamma(p)$ of all global sections. This map
is injective iff the intersection of the congruence $c_x, x\in X,$ is
the equality relation. But there is no general criterion for 
surjectivity. {\sc Surjectivity is the real challenge}.

The condition that the sets $e(a,b)$ are open can be enforced by
choosing on $X$ the coarsest topology such that all the $e(a,b)$ are
open. For example the co-Zariski topology on the prime ideals of a
commutative ring (commutative $\ell$-group, commutative $\ell$-ring)
is such a topology.

2) In case that the equalizers $e(a,b)$ are not always open, one can
enforce this condition by `localizing', that is, replacing $c_x$ by 
$\widetilde c_x=\bigcup_U\bigcap_{y\in U}c_y$, where $U$ ranges over
all neighborhoods of $x$. 

This procedure is applied when considering the prime ideals    of a
commutative ring (commutative $\ell$-group, commutative $\ell$-ring)
with the Zariski topology.  \\

{\bf Zariski or Co-Zariski topology}\\

As far as I know, until 1975, only the Zariski topology was used on the
base spaces for the sheaf representations. \\

Why the Zariski topology? Firstly, the Zariski topology was used in 
algebraic topology. Secondly: If we start with the ring $C(X)$ of
continuous real-valued functions on a compact Hausdorff space $X$, then
the topology of $X$ is reproduced by the Zariski topology on the set
of all maximal ideals. It does
not matter whether one looks at $C(X)$ as a ring or as a
lattice-ordered group. See, for example:
\begin{verse}
 R. Bkouche, Couples spectraux et faisceaux associ\'es. Applications aux
anneaux de fonctions, Bull. Soc. Math. France 98 (1970), 253 - 295.
\end{verse}
Another example is the ring direct sum $\sum_{i\in I}K_i$ of fields
$K_i$. The Zariski topology on the maximal ideals reproduces the
discrete topology on $I$, the co-Zariski topology is the cofinite
topology on $I$. \\

To the best of my knowledge, the use of the co-Zariski topology for sheaf
representations first occurs in a paper by J. F. Kennison in 1976:
\begin{verse}
J. F. Kennison, Integral domain type representations in sheaves and
other topoi. Math. Z.  151  (1976), 35-56. \\
J. F. Kennison and C. S. Ledbetter, Sheaf representations and the
Dedekind reals.  Applications of Sheaves
(Proc. Res. Sympos. Appl. Sheaf Theory to Logic, Algebra and Anal.,
Univ. Durham, Durham, 1977),  pp. 500-513, Lecture Notes in Math.,
753, Springer, Berlin-New York, 1979. 
\end{verse}
In the case of Boolean algebras of factor congruences the Zariski and
the co-Zariski topology agree on the space of maximal ideals. So there
is no difference between the two. But in other cases the Zariski
topology seems to be the natural one. I would be interesting to hear
arguments that for some purposes the co-Zariski topology is preferable. \\

{\bf Kennison's approach}\\

Kennison has an approach that is different from others and I want to
comment on it, as I understand it. (The reason is that it sheds light
on the role of the co-Zariski topology for sheaf representations.) He
considers the class $\mathcal 
A$ of all $\Omega$-algebras or an equationally defined subclass thereof.
In the class $\mathcal A$ he chooses a subclass $\Sigma$ that is
closed for subalgebras and ultraproducts. \\

Kennison asks the question: \emph{Which
algebras in $\mathcal A$ allow a representation as algebras of all
global section of a sheaf of $\Sigma$-algebras?} 

His answer is: \emph{Precisely those algebras that are in the limit
  closure of $\Sigma$.}\\

Actually I do not know whether this statement holds in general, since
Kennison suppose that the algebras in $\mathcal A$ have a group
operation among their operations. He determines this limit closure for
some examples.

1) $\mathcal A$ is the class of all unital commutative rings $R$, $\Sigma$
the subclass of integral domains. Then $R$ is representable as the
ring of all global sections of a sheaf of integral domains if and only
if $R$ has no nilpotent element AND satisfies a sequence of axioms (D$_n$).  
As a consequence, not every unital commutative ring without nilpotent
element is isomorphic to the ring of global sections of the sheaf
constructed canonically (as above) over the set of prime ideals with
the co-Zariski topology.

2) $\mathcal A$ is the class of all unital commutative f-rings $R$, $\Sigma$
the subclass of integral domains. Then $R$ is representable as the
ring of all global sections of a sheaf of integral domains if and only
if $R$ has no nilpotent element AND satisfies the axioms (D$_1$).  

3) In those early papers Kennison does not apply his general result to
the class $\mathcal A$ of unital commutative lattice ordered groups or
MV-algebras $R$ with the subclass $\Sigma$ of totally ordered ones. I am
sure that in these cases the limit closure of $\Sigma$ is all of
$\mathcal A$, which explains that in these cases one has 
representations of $R$ by the global sections of a sheaf of totally
ordered groups (MV-algebras) over the spectrum with the co-Zariski
topology. \\

{\bf Applications of sheaf representations}\\

Sheaf representations over Boolean spaces, also being rephrased as
Boolean products have been used extensively in {\bf model theory}:
\begin{verse}
A. Macintyre, Model-completeness for sheaves of structures, 
Fundamenta Mathe\-maticae 81 (1973), 73-85.\\
S. Comer, Elementary properties of structures of sections,
Bol.~Soc.~Mat. Mexicana 19 (1974), 78-85.\\
L. Lipshitz, D. Saracino, The model companion of the theory of
commutative rings without nilpotent elements,
Proc.~Amer.~Math.~Soc.~38 (1973), 381-387.\\
A. B. Carson, The model completion of the theory of commutative
regular rings, J. Algebra 27 (1973), 136-146.\\
L. van den DRIES, Artin-Schreier theory for commutative regular rings, 
Annals of Mathematical Logic 12 (1977), 113-150.\\
S. Burris, Decidability and Boolean representations, Memoirs AMS 246
(1981).\\
S. Burris, H. Werner, Sheaf constructions and their elementary
properties, Trans. Amer. Math. Soc.  248 (1979), 269-309.
\end{verse}
In later chapters of the book by Burris and Sankappanavar one can find
a well presented 
summary of model theoretic results using sheaf representations over
Boolean spaces (equivalently Boolean powers):
\begin{verse}
S. Burris, H.P. Sankappanavar, A course in universal algebra,
Springer-Verlag, 1981.
\end{verse}
Sheaf representations over Boolean spaces, equivalently Boolean
products, are particularly apt to
applying model theoretic methods. Weispfennning has tried to extend
the realm of applicability of these methods to more general sheaf
representations replacing 'Boolean' by 'distributive lattice'. But the
generalization was not too successful:
\begin{verse}
V. Weispfenning, Model theory of lattice products,
Habilitationsschrift, Universit\"at Heidelberg (235 pp.). \\
---, Lattice products, in Logic Colloquium 1978, Mons, North-Holland,
Amsterdam, pp. 423-426. 
\end{verse}
These applications in model theory are the only systematic use of
sheaf representations that I know of. Occasionally sheaf
representations have been used for exhibiting Examples with specific
properties. Already Pierce constructs \emph{modules with almost any
pathological properties that can be imagined} in the Memoir cited
above. Certain `hulls' of given algebraic structures have been
constructed using sheaf representations: 
\begin{verse}
Kist, J.: Compact spaces of minimal prime ideals. Math. Z. 111,
151-158 (1969). \\
K. Keimel, Baer extensions of rings and Stone extensions of semigroups
Semigroup Forum 2 (1971), 55-63.\\
D. H. Adams, Semigroups with not non-zero nilpotents, Math. Z. 123
(1971), 168-176.\\
B. A. Davey, m-Stone lattices, Can. J. Math., Vol. XXIV, No. 6, (1972),
1027-1032.\\ 
W. Rump, Y.~C.~Yang, Lateral completion and structure sheaf of an
archimedean l-group, J. Pure and Appl. Algebra 213 (2009), 136-143.
\end{verse}
In his paper mentioned in a previous section, B. A. Davey has given a general
explanation for the constructions given in the three papers just
mentioned. \\

Apart from these applications sheaves are mainly used for visualizing
given algebraic structures as algebras of continuous functions with
values in hopefully simpler algebras. For this reason I am not sure
whether the revival of sheaf representations in the last years is
justified. What is a representation good for, if one is not using it
as tool for obtaining new results?\\

{\bf Sheaves with values in other categories}\\

In the above I was concentrating on sheaves with values in categories
of algebras. The stalks are algebraic structures with the discrete
topology. \\

But already in early times one was considering sheaves with values in
topological spaces, uniform spaces, Banach spaces, $C^*$-algebras,
etc. One can see this already from Hofmann's SURVEY mentioned above.  
These efforts were motivated from functional analysis. Instead of the
notion of a sheaf, often the terminology `bundles of Banach spaces'
is used, although these `bundles' need not be locally trivial.\\

%
%
%
%
%
%

%
%

{\bf Questions}\\

I would like to invite every reader of these Notes to communicate to
me\footnote{keimel@mathematik.tu-darmstadt.de} any information
concerning the following questions: 
\begin{itemize}
\item Which aspects should be modified in the above?
\item Which aspects should be added to the above considerations?
\item Where have sheaf representations been used or applied to obtain results
  that do not concern the sheaf representation itself?
\item What is a justification for the use of the co-Zariski topology for sheaf
  representations?
\item Which bibliographical items  not mentioned above have
  influenced the developments in a significant way?
\item Which sheaf representation theorem published during the
  past 20 years has been used for purposes outside of sheaf
  representations?
\end{itemize}

Please keep in mind that the above reflects my personal experience and
knowledge. There must be omissions. I have not touched to questions
concerning intuitionistic 
logic, toposes etc. But I am willing to include any helpful comment.
The bibliographical items have not been selected systematically.
\end{document}